\documentclass[10pt]{article}
\usepackage{amssymb}
\usepackage{latexsym}
\usepackage{natbib}
\usepackage{amsmath}
\usepackage[all]{xy}

\usepackage[dvips]{graphicx}
\usepackage{float}
\usepackage{labelfig}
\usepackage{epsfig,psfig}

\newcommand{\B}{\Bbb}
\newcommand{\OO}{\mbox{$\mathcal{O}$}}

\long\def\symbolfootnote[#1]#2{\begingroup
\def\thefootnote{\fnsymbol{footnote}}\footnote[#1]{#2}\endgroup}

\begin{document}

\begin{center}
\textbf{\begin{large}ON THE DISTANCE BETWEEN SEIFERT SURFACES\end{large}\\~\\}

MAKOTO SAKUMA and KENNETH J. SHACKLETON\\~\\


\textit{Dedicated to Professor Takao Matumoto}

\textit{on the occasion of his sixtieth birthday\\~\\~\\}
\end{center}

\begin{center}\textbf{Abstract}\end{center}

\indent For a knot $K$ in $\B{S}^{3}$, Kakimizu introduced a simplicial complex whose vertices are all the isotopy classes of minimal genus spanning surfaces for $K$. The first purpose of this paper is to prove the $1$-skeleton of this complex has diameter bounded by a function quadratic in knot genus, whenever $K$ is atoroidal. The second purpose of this paper is to prove the intersection number of two minimal genus spanning surfaces for $K$ is also bounded by a function quadratic in knot genus, whenever $K$ is atoroidal. As one application, we prove the simple connectivity of Kakimizu's complex among all atoroidal genus $1$ knots.\\

\noindent K\begin{small}EYWORDS\end{small}: knot; atoroidal knot; Seifert surface; knot genus. \symbolfootnote[0]{2000 Mathematics Subject Classification: 57M25, 05C12}\\


\section{Introduction}

Let $K$ be a knot in the $3$-sphere $\B{S}^{3}$. A {\em Seifert surface} for $K$ is a compact, connected and orientable surface in $\B{S}^{3}$ whose boundary is precisely $K$. Fix a regular neighbourhood $N(K)$ for the knot $K$, and denote by $E(K)$, or just $E$, its exterior $\B{S}^{3} - {\rm int} N(K)$. We say that $E$, or $K$, is {\em atoroidal} if every incompressible torus in $E$ is boundary parallel. We shall say that a properly embedded subsurface of $E$ is a {\em spanning surface for $K$} if it is contained in some Seifert surface for $K$. For any spanning surface or Seifert surface $S$, we denote its ambient isotopy class by $[S]$. Throughout this paper, we shall assume, unless otherwise stated, that any given pair of Seifert surfaces or spanning surfaces intersects transversely.

To the knot $K$ there is an associated graph $\mathcal{G}(K)$, constructed as follows. The vertex set comprises the isotopy classes of minimal genus spanning surfaces for $K$, and two distinct vertices are connected by an edge if and only if they can be represented by a pair of disjoint spanning surfaces in $E$. It is a result of Scharlemann-Thompson (Proposition 5 from [13]) that $\mathcal{G}(K)$ is connected. As it happens, their main argument implies $d(\sigma, \sigma') \leq \iota(\sigma, \sigma') + 1$, where $\sigma$ and $\sigma'$ are two isotopy classes of spanning surfaces, $d$ is the path-metric on $\mathcal{G}(K)$ assigning length $1$ to each edge, and $\iota(\sigma, \sigma')$ denotes the least number of components $|S \cap S'|$ among all pairs $S \in \sigma, S' \in \sigma'$ of spanning surfaces intersecting transversely. We refer to $\iota(\sigma, \sigma')$ as the {\em intersection number of $\sigma, \sigma'$}. It is a result of Jaco-Sedgwick (see Oertel [9]) that $\mathcal{G}(K)$ is a finite graph whenever $K$ has genus at least $2$ and is atoroidal.  A result of Wilson [20] implies, among other things, that $\mathcal{G}(K)$ is a finite graph whenever $K$ is atoroidal even if it has genus $1$. However, it is a recent result of Tsutsumi [18] that there exists an infinite sequence of atoroidal knots of common genus, in each genus at least $2$, with an increasing number of isotopy classes of minimal genus spanning surfaces. In particular, the number of vertices of $\mathcal{G}(K)$ among all atoroidal knots $K$ of any common genus, at least $2$, is unbounded.

Our first result offers a uniform bound on the diameter of $\mathcal{G}(K)$, quadratic in knot genus $g(K)$.

\newtheorem{1.1}{Theorem}[section]
\begin{1.1}
Suppose that $K$ is an atoroidal knot in $\B{S}^{3}$. Then, $\mathcal{G}(K)$ has diameter at most $2g(K)(3g(K) - 2)+1$.
\label{first}
\end{1.1}

The assumption that $K$ be atoroidal is necessary here, for Theorem B of [5] asserts that the graph $\mathcal{G}(K)$ is a bi-infinite line for many composite knots $K$. Such knots are toroidal, and such graphs are unbounded.

One can go further and place a uniform quadratic bound on the intersection number of any pair of minimal genus Seifert surfaces.

\newtheorem{1.2}[1.1]{Theorem}
\begin{1.2}
Suppose that $K$ is an atoroidal knot in $\B{S}^{3}$. Then, $\iota([S], [S'])$ is at most $2(3g(K) - 2)^{2}$ for any two minimal genus Seifert surfaces $S$ and $S'$ for $K$.
\label{third}
\end{1.2}

We note Theorem \ref{third} also offers a second uniform bound on the diameter of $\mathcal{G}(K)$, though less desirable than that already given in Theorem \ref{first}. A useful consequence of Theorem \ref{third} is the following.

\newtheorem{1.1and1.2}[1.1]{Corollary}
\begin{1.1and1.2}
Suppose that $K$ is an atoroidal genus $1$ knot in $\Bbb{S}^{3}$. Then, $\mathcal{G}(K)$ has diameter at most $2$.
\label{atmost2}
\end{1.1and1.2}

The better bound on diameter offered by Corollary \ref{atmost2} for all atoroidal genus $1$ knots is in fact sharp.

\newtheorem{1.3}[1.1]{Proposition}
\begin{1.3}
There exists an atoroidal genus $1$ knot $K$ in $\B{S}^{3}$ such that $\mathcal{G}(K)$ has diameter $2$.
\label{sharp}
\end{1.3}

Recall $\mathcal{G}(K)$ is also the $1$-skeleton of the corresponding simplicial complex $\mathcal{MS}(K)$ constructed by Kakimizu [5], where a set of pairwise distinct vertices spans a simplex if and only if they may simultaneously be realised disjointly in $E$. That this complex is flag, so that any inclusion of the boundary of any given simplex extends over the whole simplex, is Proposition 4.9 from [12]. Applications of this complex are found in classifying the incompressible Seifert surfaces of prime knots of at most 10 crossings (see [6]), using a method that enhances that of Kobayashi's [8]. In [12], the first author proves the contractibility of $\mathcal{MS}(K)$ when $K$ is a special aborescent link. In his joint paper with Hirasawa [3], contractibility when $K$ is a prime, special, alternating link is announced. Together, this partially verifies a challenging conjecture of Kakimizu's [4], asserting among other things that $\mathcal{MS}(K)$ is always contractible. A full statement is given as Conjecture 0.2 in [12].

Applying Theorem \ref{third} and Corollary \ref{atmost2}, we will prove the following.

\newtheorem{1.4}[1.1]{Proposition}
\begin{1.4}
Suppose that $K$ is an atoroidal knot of genus $1$. Then, the simplicial complex $\mathcal{MS}(K)$ is simply connected.
\label{simpleconn}
\end{1.4}

We remark the dimension of $\mathcal{MS}(K)$ is at most $6$ whenever the knot is atoroidal and has genus $1$, as follows from Tsutsumi [17].

The structure of this paper is as follows. In Section 2 we shall recall Kakimizu's characterisation of the metric on $\mathcal{G}(K)$. In Section 3 we use Section 2 to prove Theorem \ref{first}. In Section 4 we prove Theorem \ref{third}, and then Corollary \ref{atmost2}. In Section 5 we prove Proposition \ref{simpleconn}. In each of Section 3, Section 4 and Section 5, we will rely on Theorem 3.1 of Fenley [1] which, among other things, rules out the existence of an accidental peripheral in any minimal genus Seifert surface whenever the knot is atoroidal. To recall, an {\em accidental peripheral} on a surface $S$ in $\B{S}^{3}$ is a simple loop essential on $S$ which is homotopic, in $E$, to a loop on $\partial E$. In Section 6 we prove Proposition \ref{sharp}. In Section 7 we investigate the growth in diameter of $\mathcal{G}(K)$ with knot genus $g(K)$, proving the following.

\newtheorem{1.5}[1.1]{Proposition}
\begin{1.5}
For every positive integer $g$, there is an atoroidal knot $K$ of genus $g$ such that the diameter of $\mathcal{G}(K)$ is equal to $2g-1$.\label{seventh}
\end{1.5}

It seems appropriate to close the introduction by posing the following open question.

\newtheorem{1.6}[1.1]{Question}
\begin{1.6}
Considering the quadratic upper bound on diameter offered by Theorem \ref{first}, can this be improved to a linear function of knot genus?
\end{1.6}

\noindent R\begin{small}EMARK:\end{small} A few months after making the first version of this paper publicly available, the authors learned in [10] that Roberto Pelayo has independently since found a version of Theorem \ref{first} as part of his PhD thesis, in preparation under the supervision of Danny Calegari.

\section{A characterisation of distance}

Let us recall Kakimizu's characterisation of the metric on $\mathcal{G}(K)$ before giving a proof to Theorem \ref{first}. For a knot $K$ in $\B{S}^{3}$ let $E$ denote its exterior and consider the infinite cyclic cover $\phi : \widetilde{E} \rightarrow E$, denoting by $\tau$ a generator for the deck transformation group. Let $S$ be any minimal genus spanning surface for $K$ and denote by $E_{0}$ the closure in $\widetilde{E}$ of any $\phi$-lift of the complement $E-S$. Set $E_{j} = \tau^{j}(E_{0})$ and $S_{j} = E_{j-1} \cap E_{j}$ for each integer $j$.

For a second such spanning surface $S'$ we may similarly form $E_{0}'$, the closure in $\widetilde{E}$ of any lift of $E-S'$ via $\phi$, and then denote by $E'_{j}$ the translate $\tau^{j}(E'_{0})$ for each integer $j$. Setting $m_{-} = {\rm {min}}\{k \in \B{Z} : E_{0} \cap E_{k}' \neq \emptyset\}$ and $m_{+} = {\rm {max}}\{k \in \B{Z} : E_{0} \cap E_{k}' \neq \emptyset\}$, we can then define $d_{*}(S, S') = m_{+} - m_{-}$. Finally, for any two vertices $\sigma, \sigma'$ of $\mathcal{G}(K)$ we define $d_{*}(\sigma, \sigma') = {\rm {min}} \{d_{*}(S, S') : S \in \sigma, S' \in \sigma'\}$. The following statement combines two key results due to Kakimizu, Proposition 1.4 and Proposition 3.2(2) of [5].

\newtheorem{2.1}[1.1]{Propostion}
\begin{2.1}{$[5]$}
The function $d_{*}$ is a metric on the vertex set of $\mathcal{G}(K)$. Moreover, for any two vertices $\sigma$ and $\sigma'$ of $\mathcal{G}(K)$, we have $d(\sigma, \sigma') = d_{*}(\sigma, \sigma')$.
\label{forth}
\end{2.1}

\section{Proof of Theorem \ref{first}}

We begin by proving two statements regarding the intersection of a pair of spanning surfaces.

\newtheorem{3.1}[1.1]{Lemma}
\begin{3.1}
Suppose that $S$ and $S'$ are two distinct minimal genus spanning surfaces for the atoroidal knot $K$. Then, $S'$ is ambient isotopic to a third minimal genus spanning surface $S''$ such that $S \cap S''$ is a disjoint union of loops and such that $d_{*}(S, S'') \leq d_{*}(S, S')$.
\label{fifthnewnew}
\end{3.1}

\noindent Proof. Suppose $S \cap S'$ contains an arc component, that is $\partial S$ and $\partial S'$ intersect. Then, $\partial S$ and $\partial S'$ bound a bigon on $\partial E$, because they are isotopic on $\partial E$. Using this bigon, we deduce $S'$ is isotopic to a second spanning surface $S''$ such that $|\partial S \cap \partial S''| \leq |\partial S \cap \partial S'| - 2$ while $d_{*}(S, S'') \leq d_{*}(S, S')$. A proof can now be completed by induction. $\square$

\newtheorem{3.2}[1.1]{Lemma}
\begin{3.2}
Suppose that $S$ and $S'$ are two distinct minimal genus spanning surfaces for the atoroidal knot $K$ intersecting only in loops. Then, $S'$ is ambient isotopic to a third minimal genus spanning surface $S''$ such that $S \cap S''$ is a disjoint union of loops, essential both on $S$ and on $S''$, and such that $d_{*}(S, S'') \leq d_{*}(S, S')$.\label{fifthnew}
\end{3.2}

\noindent Proof. Suppose for contradiction some component $\alpha$ of $S \cap S'$ is inessential on $S$, that is $\alpha$ is null-homotopic on $S$ or is boundary parallel on $S$. If $\alpha$ is null-homotopic on $S$ then, by the incompressibility of $S'$, it must also be null-homotopic on $S'$. The irreducibility of $E$ then allows us to reduce the cardinality $|S \cap S'|$ by an isotopy of $S$ or of $S'$ without increasing $d_{*}(S, S')$. We may thus assume $\alpha$ is boundary parallel on $S$. By Theorem 3.1 of [1], $\alpha$ must also be boundary parallel on $S'$.

There exist two Seifert surfaces extending $S$ and $S'$, respectively, whose intersection is precisely $(S \cap S') \cup K$. Out of convenience, we shall proceed by using $S$ and $S'$ to respectively denote such a pair of Seifert surfaces. Replacing $\alpha$ with a second component of $S \cap S'$, boundary parallel on $S$, we may assume $\alpha$ borders an annulus $A \subset S$ such that $A \cap S' = \partial A$, the union $\alpha \cup K$. Let $A' \subset S'$ be the annulus bounded by $\partial A$. The union $A \cup A'$, denoted by $T$, is an embedded torus in $\B{S}^{3}$ such that $K \subset T$.

Let $V$ be a solid torus in $\B{S}^{3}$ bounded by $T$. Since $K$ is a non-trivial and atoroidal knot, $K$ is isotopic to the core of $V$ and hence the pair $(V, A')$ can be given a product structure $(A' \times [0, 1], A' \times \{0\})$. We deduce $S'$ is isotopic to a second surface $S''$ such that $|\partial S \cap \partial S''| < |\partial S \cap \partial S'|$ while $d_{*}(S, S'') \leq d_{*}(S, S')$. A proof can now be completed by induction. $\square$\\

\noindent R\begin{small}EMARK\end{small}. We claim that in fact $S' \cap V = A'$, so that the above isotopy of $S'$ can be chosen to fix pointwise the complement of a regular neighbourhood of $A'$ in $S'$. To see this, let us argue by supposing otherwise. Let $S_{*}'$ denote the open surface $S'-A'$. The intersection of $S_{*}'$ with $V$ is then by assumption non-empty and, as $S'_{*}$ is connected and as $S'_{*}\cap T=(S'-A') \cap {\rm int} A=\emptyset$, so $S_{*}'$ is entirely contained in the interior of $V$. However, the closure of $S_{*}'$ is a surface in $V$ whose only boundary component is $\alpha$. That is, $[\alpha]$ is trivial in $H_{1}(V, \B{Z})$ despite $\alpha$ being a longitude for $V$. We have a contradiction, and we deduce $S' \cap V = A'$.\\

We shall henceforth denote by $\eta$ the function $3g - 2$ on the set of all knots, noting $\eta(K)$ is the size of any maximal collection of pairwise disjoint and non-isotopic essential simple loops on any minimal genus spanning or Seifert surface for $K$. Since an estimate on diameter is easily found for a trivial knot, we shall also assume the genus of $K$ is positive.

Let $\sigma$ and $\sigma'$ be any two vertices of $\mathcal{G}(K)$, and take representatives $S \in \sigma$ and $S' \in \sigma'$ so that $d_{*}(S, S')$ is minimal and so that $S$ and $S'$ intersect transversely in a disjoint union of essential simple loops, as per Lemma \ref{fifthnewnew} and Lemma \ref{fifthnew}.

Suppose for contradiction $d(\sigma, \sigma') \geq 2g(K)\eta(K) + 2$. According to Proposition \ref{forth}, we also have $d_{*}(S, S') \geq 2g(K)\eta(K) + 2$. Let $Y_{j}$ denote the surface $E_{j} \cap S'_{0}$ for each integer $j$ where, perhaps after reindexing, we may assume $Y_{j}$ is non-empty if and only if $1 \leq j \leq d_{*}(S, S')$. It should be noted $\chi(Y_{j}) \leq 0$ for each such index $j$, so that $|\chi(S_{0}')| = \sum_{j=1}^{d_{*}(S, S')} |\chi(Y_{j})|$. Here, $\chi(S)$ denotes the Euler characteristic of a surface $S$.\\

\noindent C\begin{small}LAIM\end{small}. \textit{For the finite sequence $\{Y_{j} : 2 \leq j \leq d_{*}(S, S') - 1\}$ of non-empty surfaces, there are $\eta(K)$ consecutive indices whose corresponding surfaces each have zero Euler characteristic. That is, there exists a natural number $r$, with $1 \leq r \leq d_{*}(S, S') - \eta(K) - 1$, such that $\chi(Y_{j}) = 0$ for each $j$ with $r + 1 \leq j \leq r + \eta(K)$.\\}

\noindent Proof.  Set $||K||$ equal to $2g(K) - 1$, the Thurston norm [15] of either generator of $H_{2}(E, \partial E)$. In the identity $||K|| = |\chi(S_{0}')| = \sum_{j=1}^{d_{*}(S, S')} |\chi(Y_{j})|$, only at most $||K||$ of the summands $|\chi(Y_{j})|$ can be non-zero. The proof will now be completed by a pigeonhole-type argument, in the following manner.

We denote by $w$ the string $(w_{2}, \ldots, w_{d_{*}(S, S') - 1})$ of binary digits, where $w_{j}$ is defined equal to $0$ if $\chi(Y_{j})$ is $0$ or otherwise $1$ for each of the indices $j \in \{2, \ldots, d_{*}(S, S') - 1\}$. If it should happen that for any $\eta(K)$ consecutive binary digits $w_{j}$ at least one is always non-zero, we would then have the  estimate $|w| \leq (||K|| + 1)\eta(K) - 1$ on the length $|w|$ of $w$. This follows from the fact that only at most $||K||$ of the binary digits $w_{j}$ can be non-zero. We can now find an upper bound for $d(\sigma, \sigma')$ as follows:

\begin{eqnarray*}
d(\sigma, \sigma') &=& d_{*}(S, S') \\
 &=& |w| + 2\\
 &\leq& (||K|| + 1)\eta(K) + 1 \\
 &=& 2g(K)\eta(K) + 1. \\
\end{eqnarray*}

\noindent According to our standing assumption on $d(\sigma, \sigma')$, this is absurd. We deduce that there exist $\eta(K)$ consecutive zeros $w_{r + 1}, \ldots, w_{r + \eta(K)}$, thus proving the claim. $\square$\\

After shifting the indexing $E_{j}$ by $r$, we have $|\chi(Y_{j})| = 0$ for $1 \leq j \leq \eta(K)$, and the set $Y_{1} \cup \cdots \cup Y_{\eta(K)}$ is both non-empty and a union of pairwise disjoint annuli. Note that both $Y_{0}$ and $Y_{\eta(K) + 1}$ are necessarily non-empty.  In particular, both $S_{1} \cap S_{0}'$ and $S_{\eta(K) +1} \cap S_{0}'$ are non-empty. There thus exists a subannulus $A$ of $Y_{1} \cup \cdots \cup Y_{\eta(K)}$ ending on $S_{1}$ and on $S_{\eta(K)+1}$.\\

\noindent C\begin{small}LAIM\end{small}. \textit{There exist natural numbers $p$ and $q$, with $1 \leq p < q \leq \eta(K) + 1$, for which there is a component $\alpha$ of $A \cap S_{p}$ and a component $\beta$ of $A \cap S_{q}$ such that $\phi(\alpha)$ and $\phi(\beta)$ are isotopic loops on $S$.\\}

\noindent Proof. For each $1 \leq j \leq \eta(K) + 1$, $\phi(A \cap S_{j})$ is a non-empty collection of essential, pairwise disjoint and pairwise non-isotopic loops on $S$, and the two sets $\phi(A \cap S_{i})$ and $\phi(A \cap S_{j})$ are disjoint for distinct $i$ and $j$. Since any collection of pairwise disjoint and non-isotopic essential simple loops on $S$ has size at most $\eta(K)$, we deduce the claim. $\square$\\

Let $\mathcal{A}$ be the family of all those subannuli of $A$ bounded by any pair of loops found in the previous claim. Then, $\mathcal{A}$ is non-empty and we can choose $A' \in \mathcal{A}$ minimal subject to inclusion. The annulus $A'' \subset S$ bounded by $\phi(\partial A')$ has interior disjoint from $\phi(A')$ and so the union $\phi(A') \cup A''$, denoted $T$, is an embedded torus in $E$.\\

\noindent C\begin{small}LAIM\end{small}. \textit{$T$ is incompressible in $E$.}\\

\noindent Proof.  We shall check the inclusion $e : T \rightarrow E$ induces an injection $e_{*} : \pi_{1}(T) \rightarrow \pi_{1}(E)$ on fundamental groups. Let $\alpha$ be either component of $\partial A'$, and let $p < q$ be such that $\partial A' \subset S_{p} \cap S_{q}$. Let the simple loop $\gamma \subset T$ be the union of an arc in $\phi(A')$ and an arc in $A''$. Observe that $\pi_{1}(T)$ is generated by $\phi(\alpha)$ and $\gamma$, and that the image of $\phi(\alpha)$ and the image of $e_{*}(\gamma)$ in $H_{1}(E) \cong \B{Z}$ are $0$ and $q - p$, respectively. Since $q - p$ is non-zero, it follows $Ker(e_{*})$ is contained in the group $\langle \phi(\alpha) \rangle$. However, $\phi(\alpha)$ is an essential loop on the incompressible surface $S$ in $E$. Hence $Ker(e_{*})$ is trivial, and $T$ is incompressible in $E$. $\square$\\

\noindent C\begin{small}LAIM\end{small}. \textit{$T$ is essential in $E$.\\}

\noindent Proof. The loop $\phi(\alpha) \subset T$ is essential on $S$. It follows from Theorem 3.1 [1] that $\phi(\alpha)$ is not isotopic in $E(K)$ to a simple loop on $\partial E$. Hence $T$ can not be boundary parallel in $E$. $\square$\\

We now have a contradiction, for $K$ is atoroidal, and we deduce $2g(K)\eta(K) + 1$ is an upper bound for the diameter of $\mathcal{G}(K)$. This completes a proof of Theorem \ref{first}.

\section{Proof of Theorem \ref{third}}

We shall once more denote by $\eta(K)$ the number $3g(K) - 2$. An argument similar to that found in the proof of both Lemma \ref{fifthnewnew} and Lemma \ref{fifthnew} permits us to represent any given pair of vertices of $\mathcal{G}(K)$ by a pair of spanning surfaces for $K$ intersecting transversely and minimally, up to isotopy, in loops essential on both surfaces.

Let $S$ and $S'$ be a pair of such spanning surfaces. Suppose for contradiction that $|S \cap S'| \geq 2\eta(K)^{2} + 1$. Then, there exist two distinct annuli $A \subset S$ and $A' \subset S'$ such that $|A \cap A'| = 3$ and $\partial A \cup \partial A' \subset A \cap A'$. To see this, consider an $l \times m$ array, with $1 \leq l, m \leq \eta(K)$, whose entries are non-negative integers summing to $2\eta(K)^{2} + 1$. It is not so hard to see that at least one of these entries must be at least $3$. We may further assume $A$ is minimal subject to inclusion, so that no component of $A \cap (S' - A')$ is isotopic on $S'$ to the core of $A'$.

The union $A \cup A'$ always separates the $3$-sphere into three components, whose closures we denote by $X_{1}, X_{2}, X_{3}$ and indexed so that $K \subset X_{3}$. Note it can happen that one of the $X_{i}$ fails to be a manifold, in which case its frontier, ${\rm fr} X_{i}$, is homeomorphic to an immersed torus whose singular set is the simple loop ${\rm int} A \cap {\rm int} A'$.\\

\noindent C\begin{small}LAIM\end{small}. \textit{${\rm int} X_{3}$ is not an open solid torus.\\}

\noindent Proof. Suppose for contradiction ${\rm int} X_{3}$ is an open solid torus. As the knot $K$ is atoroidal, so either $K$ is contained in a compact $3$-ball inside ${\rm int} X_{3}$ or $K$ is a core of ${\rm int} X_{3}$, and we rule out both cases separately as follows.\\

\noindent C\begin{small}ASE\end{small} I. \textit{$K$ is contained in a compact $3$-ball $B \subset {\rm int} X_{3}$.}  Since $S$ is connected and since $S \cap {\rm fr} X_{3}$ contains a simple loop essential on $S$, so $S \cap \partial B$ also contains a simple loop essential on $S$. Thus, there exists a disc $D$ disjoint from the knot $K$ whose boundary $\partial D$ is a non-trivial simple loop on $S$. However, $S$ is incompressible and we therefore have a contradiction.\\

\noindent C\begin{small}ASE\end{small} II. \textit{$K$ is a core of ${\rm int} X_{3}$.} Let $F$ denote the open surface $S \cap ({\rm int} X_{3} - N(K))$. Then, $F$ is necessarily a non-empty disjoint union of open annuli, for the inclusion of $F$ in $E$ descends to a monomorphism on fundamental groups that factors through the abelian group $\pi_{1}({\rm int} X_{3} - N(K)) \cong \Bbb{Z} \oplus \Bbb{Z}$. Thus, $S$  contains an annulus with one boundary component equal to $K$ and the other a component of $S \cap S'$. That is, $S$ and $S'$ intersect in at least one simple loop peripheral on $S$. However, this is contrary to the standing assumption that $S$ and $S'$ intersect only in loops essential on $S$ (and on $S'$). $\square$\\

To complete a proof of Theorem \ref{third} it suffices to rule out the following two mutually exclusive cases. These correspond to the two distinct ways in which $A$ and $A'$ can intersect one another.\\

\noindent C\begin{small}ASE\end{small} I. \textit{${\rm fr} X_{1}, {\rm fr} X_{2}, {\rm fr} X_{3}$ are each tori.} Then, at least one of $X_{1} \cup X_{2}$ and $X_{3}$ is a solid torus and, according to the claim, it can only be $X_{1} \cup X_{2}$. It follows that both $X_{1}$ and $X_{2}$ are solid tori, and, using van Kampen's theorem, at least one, say $X_{1}$, has a product structure $(X_{1}, A'') \cong (A'' \times [0, 1], A'' \times \{0\})$, where $A'' \subseteq {\rm fr} X_{1} \cap S$ is an annulus. Thus, $S$ is ambient isotopic to a second Seifert surface intersecting $S'$ fewer than $\iota([S], [S'])$ times and this is absurd.\\

\noindent C\begin{small}ASE\end{small} II. \textit{Exactly one of ${\rm fr} X_{1}, {\rm fr} X_{2}, {\rm fr} X_{3}$ is not an embedded torus.} Then, either ${\rm fr} X_{3}$ is not an embedded torus or exactly one of ${\rm fr} X_{1}$ or ${\rm fr} X_{2}$ is not an embedded torus. We thus need only consider the following two subcases.\\

\noindent II.1. \textit{${\rm fr} X_{3}$ is an immersed singular torus.}  Then, $N(X_{1} \cup X_{2})$ is necessarily a solid torus. It follows both $X_{1}$ and $X_{2}$ are solid tori and at least one, say $X_{1}$, has a product structure $(X_{1}, A'') \cong (A'' \times [0, 1], A'' \times \{0\})$, where $A'' \subseteq {\rm fr} X_{1} \cap S$ is an annulus. Once more, we will find $S$ is ambient isotopic to a second Seifert surface intersecting $S'$ fewer than $\iota([S], [S'])$ times and this is absurd.\\

\noindent II.2. \textit{${\rm fr} X_{3}$ is an embedded torus.} Then, $X_{1} \cup X_{2}$ is necessarily a solid torus. Let $F'$ denote that component of $S' \cap (X_{1} \cup X_{2})$ containing ${\rm int} A \cap {\rm int} A'$. The minimality of $A$ implies $\partial F' - \partial A$ comprises of loops essential on $S'$, none of which is isotopic on $S'$ to the core of $A'$. Hence $\pi_{1}(F')$ is non-abelian and, as $S'$ is incompressible, the inclusion of $F'$ in $E$ descends to a monomorphism on fundamental groups. In particular, it has non-abelian image. However, this monomorphism also factors through the abelian group $\pi_{1}(X_{1} \cup X_{2}) \cong \Bbb{Z}$ and as such has abelian image, a contradiction.\\

We thus complete a proof of Theorem \ref{third}. Let us finish this section by providing a proof of Corollary \ref{atmost2}.\\

\noindent Proof of Corollary \ref{atmost2}: Let $\sigma_{1}$ and $\sigma_{2}$ be two vertices of $\mathcal{G}(K)$. Let $S_{1} \in \sigma_{1}$ and $S_{2} \in \sigma_{2}$ be a pair of representative spanning surfaces, together realising intersection number, and such that $S_{1} \cap S_{2}$ is a collection of loops, perhaps empty, essential both on $S_{1}$ and on $S_{2}$. Then, since $g(S_{1})$ and $g(S_{2})$ are both equal to $1$, so $S_{1} \cap S_{2}$ comprises only of non-separating loops, parallel on $S_{1}$ and on $S_{2}$. Applying to Theorem \ref{third}, we have $|S_{1} \cap S_{2}| \leq 2$ and it follows that each lift of $S_{2}$ intersects at most one lift of $S_{1}$. Thus, $d_{*}(S_{1}, S_{2}) \leq 2$. By Proposition \ref{forth}, we have $d(\sigma_{1}, \sigma_{2}) \leq 2$ as required. $\square$\\

\section{Proof of Proposition \ref{simpleconn}}

We shall need the following criterion for the simple connectivity of a simplicial complex whose $1$-skeleton is a metric graph of diameter at most $2$, and then a restriction on the intersection number of two genus $1$ spanning surfaces.

\newtheorem{5.1}[1.1]{Lemma}
\begin{5.1}
Suppose $\mathcal{C}$ is a simplicial complex, whose $1$-skeleton can be realised as a  metric graph of diameter at most $2$, for which every simplicial circuit of length at most $5$ is contractible. Then, $\mathcal{C}$ is simply connected.\label{criterion}
\end{5.1}

\noindent Proof. Let $\sigma_{1}, \ldots, \sigma_{n}$ be the cyclically indexed vertices of a circuit $c$ of length $n$. Since $d(\sigma_{1}, \sigma_{i}) \leq 2$ for $3 \leq i \leq n-1$, so there exists a simplicial path of length at most $2$ connecting $\sigma_{1}$ to $\sigma_{i}$ for each such $i$. It follows $c$ can be expressed as a finite sum of simplicial $3$-, $4$- and $5$-circuits. Each such circuit is contractible, by assumption, and so $c$ must also be contractible. Hence, $\mathcal{C}$ is simply connected. $\square$\\

By the proof of Corollary \ref{atmost2}, we have the following lemma.

\newtheorem{5.2}[1.1]{Lemma}
\begin{5.2}
Let $K$ be an atoroidal genus $1$ knot. Then, for any pair of vertices $\sigma$ and $\sigma'$ of $\mathcal{G}(K)$, we have $\iota(\sigma, \sigma') \in \{0, 2\}$.
\label{only2}
\end{5.2}

In proving the following lemma, we shall make use of a construction that amounts to a special case of the so-called \textit{double curve sum}, after Scharlemann-Thompson [13], and of a construction of Kakimizu's [5].

\newtheorem{5.3}[1.1]{Lemma}
\begin{5.3}
Let $K$ be an atoroidal genus $1$ knot. For any pair of vertices $\sigma$ and $\sigma'$ of $\mathcal{G}(K)$, with $d(\sigma, \sigma') = 2$, there exists a third vertex $\delta$ such that $d(\delta, \sigma) = d(\delta, \sigma') = 1$ and such that, for any fourth vertex $\mu$, if $\iota(\mu, \sigma) = \iota(\mu, \sigma') = 0$, then $\iota(\mu, \delta) = 0$.\label{tight}
\end{5.3}

\noindent Proof. By Lemma \ref{only2} and arguments given in Section 3, there exist representatives $S \in \sigma$ and $S' \in \sigma'$ such that $S \cap S'$ is a pair of loops essential both on $S$ and on $S'$. Let $P \subset S$ denote the $3$-holed sphere bordered by $\partial S$ and by $S \cap S'$, and let $A' \subset S'$ be the closed annulus bordered by $S \cap S'$. Then, $P \cup A'$ is a genus $1$ spanning surface for the knot $K$ and, after a small isotopy, is disjoint from both $S$ and $S'$. We take $\delta$ to be the isotopy class $[P \cup A']$, noting that $\iota(\delta, \sigma) = \iota(\delta, \sigma') = 0$ by construction. Since $d(\sigma, \sigma') = 2$, so $d(\delta, \sigma) = d(\delta, \sigma') = 1$.

Now suppose $\mu$ is a fourth vertex, adjacent to both $\sigma$ and $\sigma'$. We claim that $\iota(\delta, \mu) = 0$, and to prove this it suffices to prove the existence of a representative of $\mu$ simultaneously disjoint both from $S$ and from $S'$. By assumption, $\iota(\mu, \sigma) = 0$ and hence there exists a representative $R \in \mu$ disjoint from $S$. Perhaps after replacing $R$ with an isotopic surface also disjoint from $S$, we may further assume that $R$ is transverse to $S'$.

As $\iota(\mu, \sigma') = 0$, by Proposition 4.8(2) of [12] there exists a product region $V$ between $R$ and $S'$ such that $V \cap R = {\rm fr} V \cap R$ and $V \cap S' = {\rm fr} V \cap S'$. Note that, should $S \cap V$ not be empty, then $S \cap V$ is parallel in $V$ to a subsurface of $S' \cap V$, by Corollary 3.2 of [19]. It follows $S$ and $S'$ would share a removable intersection, and this is absurd. Thus, $S \cap V$ is empty. We can therefore use the region $V$ to replace $R$ with an isotopic surface $R'$ such that $R'$ and $S$ are disjoint and such that $|R' \cap S'| \leq |R \cap S'| - 1$.

Continuing inductively, we deduce $R$ is isotopic to a spanning surface simultaneously disjoint both from $S$ and from $S'$. $\square$\\

\newtheorem{5.4}[1.1]{Lemma}
\begin{5.4}
Let $K$ be an atoroidal genus $1$ knot. Suppose $\sigma, \sigma_{1}$ and $\sigma_{2}$ are three vertices of $\mathcal{G}(K)$ such that $d(\sigma, \sigma_{1}) = d(\sigma, \sigma_{2}) = 2$ and such that $d(\sigma_{1}, \sigma_{2}) = 1$. Then, there exist two vertices $\delta_{1}$ and $\delta_{2}$ of $\mathcal{G}(K)$ such that $d(\sigma, \delta_{i}) = d(\delta_{i}, \sigma_{i}) = 1$, for $i \in \{1, 2\}$, and such that $d(\delta_{1}, \delta_{2}) \leq 1$.
\label{3isok}
\end{5.4}

\noindent Proof. Let $S \in \sigma, S_{1} \in \sigma_{1}$ and $S_{2} \in \sigma_{2}$ be such that $S_{1} \cap S_{2}$ is empty and such that $S$ intersects both $S_{1}$ and $S_{2}$ transversely and in a collection of loops essential on each surface.

Let $\widetilde{E}$ denote the infinite cyclic cover of the knot exterior $E$, with covering map denoted $\phi$, and denote by $\tau$ either generator of the deck transformation group. Let $\widetilde{S}_{1,0}$ denote any lift of $S_{1}$, and let $\widetilde{S}_{1,n}$ denote the translate $\tau^{n}(\widetilde{S}_{1,0})$ for each integer $n \in \B{Z}$. We similarly introduce the notation $\widetilde{S}_{2,n}$, where $\widetilde{S}_{2,0}$ is to separate $\widetilde{S}_{1,0}$ and $\widetilde{S}_{1,1}$.

The following claim permits us to isotope $S$ so that in addition each lift of $S$ intersects only one lift of $S_{1}$ and only one lift of $S_{2}$. Recall the definition of the function $d_{*}$ from Section 2.\\

\noindent C\begin{small}LAIM\end{small}. There exists an isotopy of $S$ after which $d_{*}(S, S_{1}) = d_{*}(S, S_{2}) = 2$.\\

\noindent Proof. Suppose $d_{*}(S, S_{1}) + d_{*}(S, S_{2}) \geq 5$, and denote by $\widetilde{S}$ any lift of $S$. Then, since $d_{*}(\sigma, \sigma_{1}) = d_{*}(\sigma, \sigma_{2}) = 2$, so there exists a component $R$ of $\phi^{-1}(S_{1} \cup S_{2})$ such that $R \cap \widetilde{S}$ is not empty and for which there exists an isotopy of $S$ lifting to an isotopy of $\widetilde{S}$ after which $R \cap \widetilde{S}$ is empty. By Proposition 4.8(2) of [12], there thus exists a product region $V$ in $\widetilde{E}$ between $R$ and $\widetilde{S}$ and such that $V \cap R = {\rm fr} V \cap R$ and $V \cap \widetilde{S} = {\rm fr} V \cap \widetilde{S}$. As $\widetilde{S}$ and each component of $\phi^{-1}(S_{1} \cup S_{2})$ separate in $\widetilde{E}$, so there exists a subregion $V' \subseteq V$ such that $V' \cap (\widetilde{S} \cup \phi^{-1}(S_{1} \cup S_{2})) = {\rm fr} V'$. Let us denote by $R'$ the one component of $\phi^{-1}(S_{1} \cup S_{2})$ such that $R' \cap V'$ is not empty.

Applying Corollary 3.2 of [19] to the product region $V$, we find $R' \cap V'$ and $\widetilde{S} \cap V'$ are parallel through $V'$. Note, $V'$ is contained in a single fundamental region. Projecting $V'$ to $E$ then, we can therefore isotope $S$ so as to remove the corresponding intersection between $R'$ and $\widetilde{S}$ and without introducing any new intersections between $\phi^{-1}(S_{1} \cup S_{2})$ and $\widetilde{S}$.

That is, so long as $d_{*}(S, S_{1}) + d_{*}(S, S_{2}) \geq 5$, we can successively remove intersections between $\phi^{-1}(S_{1} \cup S_{2})$ and $\widetilde{S}$ via an isotopy of $S$. There are only finitely many such intersections to begin with, thus in finite time we construct an isotopy of $S$ after which $d_{*}(S, S_{1}) + d_{*}(S, S_{2}) \leq 4$. The statement of the claim is deduced. $\square$\\

Isotope $S$ as indicated by the claim, and denote by $\widetilde{S}$ the lift of $S$ intersecting $\widetilde{S}_{1,0}$ and $\widetilde{S}_{2,0}$. Now let $N$ denote a small regular neighbourhood of $\widetilde{S} \cup \widetilde{S}_{1,0} \cup \widetilde{S}_{2,0}$ in the infinite cyclic cover $\widetilde{E}$, so that $\tau^{j}(N)$ is disjoint from $\widetilde{S}, \widetilde{S}_{1,0}$ and $\widetilde{S}_{2,0}$ for each non-zero integer $j$. We define $\widetilde{R}_{1}$ and $\widetilde{R}_{2}$ to be the two ``outermost'' components of ${\rm fr} N$, that is $\widetilde{R}_{1}$ and $\widetilde{R}_{2}$ bound a region in $\widetilde{E}$ containing $N$ and indexed so that $\widetilde{R}_{1}$ and $\widetilde{S}_{2,0}$ are separated by $\widetilde{S}_{1,0}$. Note, $\widetilde{R}_{1}$ and $\tau^{-1}(\widetilde{R}_{2})$ are disjoint, are both $1$-holed tori, and are both contained in the fundamental region bordered by $\widetilde{S}_{1,-1}$ and $\widetilde{S}_{1,0}$. Thus, $\widetilde{R}_{1}$ and $\tau^{-1}(\widetilde{R}_{2})$ project to disjoint genus $1$ spanning surfaces, denoted respectively $R_{1}$ and $R_{2}$, both of which are disjoint from $S_{1}$ and from $S_{2}$.

Finally, we respectively define $\delta_{1}$ and $\delta_{2}$ to be the isotopy classes $[R_{1}]$ and $[R_{2}]$. This completes a proof of Lemma \ref{3isok}. $\square$\\

In view of Corollary \ref{atmost2} and Lemma \ref{criterion}, to prove the simple connectivity of $\mathcal{MS}(K)$, for an atoroidal genus $1$ knot $K$, it suffices to prove the following three claims.\\

\noindent C\begin{small}LAIM\end{small}. \textit{Every simplicial $3$-circuit in $\mathcal{MS}(K)$ is contractible.\\}

\noindent Proof. This is immediate, for $\mathcal{MS}(K)$ is a flag simplicial complex. That is, any embedding of the $1$-skeleton of any given simplex into $\mathcal{G}(K)$ is the restriction of an embedding from the whole simplex into $\mathcal{MS}(K)$. $\square$\\

\noindent C\begin{small}LAIM\end{small}. \textit{Every simplicial $4$-circuit in $\mathcal{MS}(K)$ is contractible.\\}

\noindent Proof. Suppose $\sigma_{1}, \sigma_{2}, \sigma_{3}, \sigma_{4}$ are the cyclically indexed vertices of a simplicial $4$-circuit in $\mathcal{MS}(K)$. Assuming $d(\sigma_{1}, \sigma_{3}) = 2$, by Lemma \ref{tight} there exists a vertex $\delta$ such that $d(\delta, \sigma_{1}) = d(\delta, \sigma_{3}) = 1$ and such that $\iota(\delta, \sigma_{2}) = \iota(\delta, \sigma_{4}) = 0$. We deduce $\delta$ spans an edge with $\sigma_{2}$ and with $\sigma_{4}$. Appealing to the previous claim, one may now find an appropriate compressing disc as the union of at most four $2$-simplices.

The remaining cases may be similarly treated. $\square$\\

\noindent C\begin{small}LAIM\end{small}. \textit{Every simplicial $5$-circuit in $\mathcal{MS}(K)$ is contractible.\\}

\noindent Proof. Suppose $\sigma_{1}, \sigma_{2}, \sigma_{3}, \sigma_{4}, \sigma_{5}$ are the cyclically indexed vertices of a simplicial $5$-circuit in $\mathcal{MS}(K)$. Assuming $d(\sigma_{1}, \sigma_{3}) = d(\sigma_{1}, \sigma_{4}) = 2$, by Lemma \ref{3isok} there exist vertices $\delta_{3}$ and $\delta_{4}$ of $\mathcal{G}(K)$ such that $d(\sigma_{1}, \delta_{i}) = d(\delta_{i}, \sigma_{i}) = 1$, for $i \in \{3, 4\}$, and such that $d(\delta_{3}, \delta_{4}) \leq 1$. If furthermore $\delta_{3}$ and $\delta_{4}$ are distinct, then $(\sigma_{1}, \delta_{3}, \delta_{4})$ is a $3$-circuit and $(\sigma_{1}, \sigma_{2}, \sigma_{3}, \delta_{3}), (\delta_{3}, \sigma_{3}, \sigma_{4}, \delta_{4})$, and $(\sigma_{1}, \delta_{4}, \sigma_{4}, \sigma_{5})$ are each circuits of length at most $4$. Appealing to the previous two claims respectively, one may now find an appropriate compressing disc as the union of at most four other discs.

The remaining cases may be similarly treated. $\square$\\

This completes a proof of Proposition \ref{simpleconn}.

\section{A genus 1 knot}

The purpose of this section is to prove Proposition \ref{sharp}, that is to construct an atoroidal genus $1$ knot $K$ whose graph $\mathcal{G}(K)$ has diameter $2$.

Let $V_{0}$ be a solid torus, and let $A_{1}$ and $A_{2}$ be annuli on $\partial V_{0}$ essential in $V_{0}$ such that $\partial A_{1} \cap \partial A_{2} = \partial A_{1} = \partial A_{2}$ and such that the cyclic group $\pi_{1}(V_{0})$ is not generated by the core of $A_{1}$. Note, $\partial A_{1} \cup \partial A_{2} = \partial V_{0}$. Let $V$ be a genus $2$ handlebody obtained from $V_{0}$ by attaching a $1$-handle $D^{2} \times [1, 2]$, where $D^{2} \times \{i\}$ is identified with a disc in ${\rm int} A_{i}$ for both $i \in \{1, 2\}$. By assumption, the region in $V$ bounded by $A_{1}$ and $A_{2}$ does not admit a product structure $A_{1} \times [0, 1]$. After pushing ${\rm int} A_{i}$ into ${\rm int} V$, for both $i \in \{1, 2\}$, we have a pair of annuli properly embedded in $V$.

\begin{figure}[h]
\leavevmode \SetLabels
\L(0.62*0.76) $\alpha$\\
\L(0.51*0.934) $\beta$\\
\L(0.30*0.855) $k$\\
\endSetLabels
\begin{center}
\AffixLabels{\centerline{\epsfig{file = 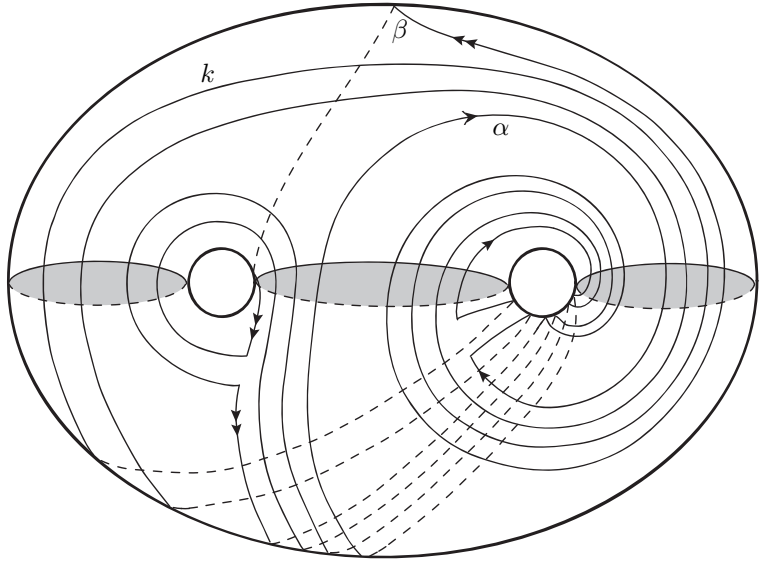, width=10cm, angle= 0}}}
\end{center}
\caption{The curve $k$ is a band sum of $\alpha$ and $\beta$, and is complicated with respect to the indicated maximal meridan system in the sense of Kobayashi.}
\label{fig:newknot}
\end{figure}

Let $\alpha, \beta$ be the two components of $\partial A_{1}$, with orientation induced by either orientation of $A_{1}$ (see Figure \ref{fig:newknot}). We can choose $A_{1}$ and a band sum $k$ of $\alpha$ and $\beta$ such that $k$ is \textit{complicated} with respect to a preferred maximal meridan system $\{D_{1}, D_{2}, D_{3}\}$ for $V$, in the sense of Kobayashi [7]. That is, if $R_{1}$ and $R_{2}$, both compact $3$-holed spheres, denote the two components of $\partial V - {\rm int} N(\partial D_{1} \cup \partial D_{2} \cup \partial D_{3})$, then:

\begin{itemize}
\item There is no bigon $B$ on $\partial V$ such that $\partial B$ is the union of a subarc of $k$ and a subarc of $\partial R_{i}$ for some $i \in \{1, 2\}$, and

\item For any two boundary components of either $3$-holed sphere $R_{i}$, there is a subarc of $k$ joining them in $R_{i}$.
\end{itemize}

Appealing to Lemma 6.1 of [7], we have the following.

\newtheorem{6.1}[1.1]{Lemma}
\begin{6.1}
$\partial V - {\rm int} N(k)$ is incompressible in $V$.
\label{vincmp}
\end{6.1}

Perhaps after an isotopy, we may assume $k$ is disjoint from $\alpha \cup \beta$. Let $\theta$ denote any graph with two vertices, connected by three edges, embedded in $\B{S}^{3}$ and whose exterior $W$ admits a complete and finite volume hyperbolic metric with totally geodesic boundary. According to Section $3.3$ of [16], one may, for instance, take $\theta$ to be the Kinoshita theta curve. Let $f : V \rightarrow N(\theta)$ be any homeomorphism, and define $K$ equal to the image $f(k)$. Let $N(K)$ be a regular neighbourhood of $K$ in $\B{S}^{3}$ such that $N(K) \cap f(V)$ and $N(K) \cap W$ are regular neighbourhoods of $K$ in $f(V)$ and $W$ respectively. Note, the exterior $E = \Bbb{S}^{3} - {\rm int} N(K)$ of $K$ is the union of $f(V) - {\rm int} N(K)$ and $W - {\rm int} N(K)$, with common subsurface $\partial f(V) - {\rm int} N(K) = \partial W - {\rm int} N(K)$. It should also be noted that $f(V) - {\rm int} N(K)$ is homeomorphic to $f(V)$ via a homeomorphism constant on $\partial f(V) - {\rm int} N(K)$, and that $W - {\rm int} N(K)$ is homeomorphic to $W$ via a homeomorphism constant on $\partial W - {\rm int} N(K)$.

\newtheorem{6.2}[1.1]{Lemma}
\begin{6.2}
The surface $\partial W - {\rm int} N(K)$ is incompressible in $W$.
\label{wincmp}
\end{6.2}

\noindent Proof. Since $K$ is an essential loop on $\partial W$ so the natural inclusion $\partial W - {\rm int} N(K) \rightarrow \partial W$ descends to an injection on fundamental groups. As $W$ admits a hyperbolic metric in which $\partial W$ is totally geodesic so $\partial W$ is incompressible in $W$, and we find the natural homomorphism $\pi_{1}(\partial W - {\rm int} N(K)) \rightarrow \pi_{1}(W)$ is also injective. It follows $\partial W - {\rm int} N(K)$ is incompressible in $W$. $\square$\\

\newtheorem{6.3}[1.1]{Lemma}
\begin{6.3}
$K$ is a non-trivial knot in $\B{S}^{3}$.
\label{nontriv}
\end{6.3}

\noindent Proof. According to Lemma \ref{vincmp}, the group $\pi_{1}(\partial V - {\rm int} N(K))$ naturally injects into $\pi_{1}(V)$. According to Lemma \ref{wincmp}, the same group $\pi_{1}(\partial V - {\rm int} N(K))$ naturally injects into $\pi_{1}(W)$. The knot group $\pi_{1}(E)$ is, by using van Kampen's theorem, therefore isomorphic to the amalgamated free product of $\pi_{1}(V)$ and $\pi_{1}(W)$ over a common subgroup isomorphic to the fundamental group of a $2$-holed torus. Hence, $\pi_{1}(E)$ is a non-abelian group, and $K$ cannot be a trivial knot. $\square$\\

\newtheorem{6.4}[1.1]{Lemma}
\begin{6.4}
The pair $(W, \partial W - {\rm int} N(K))$ does not contain an essential annulus. That is, suppose $A$ is an essential annulus properly embedded in $W$, with $\partial A \subset \partial W - {\rm int} N(K)$. Then, $A$ is parallel to an annulus in $\partial W - {\rm int} N(K)$ or to the annulus $N(K) \cap \partial W$.
\label{annuli}
\end{6.4}

\noindent Proof. The pair $(W, \partial W)$ cannot contain an essential annulus, for $W$ admits a hyperbolic metric in which $\partial W$ is totally geodesic. Let $A$ be any incompressible annulus properly embedded in $W - {\rm int} N(K)$ so that $\partial A \subset \partial W - {\rm int} N(K)$. Then, $A$ is parallel in $W$ to an annulus $A' \subset \partial W$. If $A' \cap K$ is empty then $A$ is also parallel to $A'$ in $W - {\rm int} N(K)$. If instead $A' \cap K$ is not empty, then $K \subset A'$ and, as $K$ is essential on $\partial W$, so $K$ is the core of $A'$. Thus, $A$ is parallel to the annulus $N(K) \cap \partial W$. $\square$\\

\newtheorem{6.5}[1.1]{Lemma}
\begin{6.5}
$K$ is an atoroidal knot in $\B{S}^{3}$.
\end{6.5}

\noindent Proof. Suppose $T$ is an incompressible torus in $E$. As both $f(V)$ and $W$ are atoroidal, we may assume that $T$ intersects $\partial W - {\rm int} N(K)$ only in a collection of loops essential on $\partial W$ and that each component of $T \cap f(V)$ and $T \cap W$ is an incompressible annulus in $f(V)$ and $W$, respectively.

Let $A$ be a component of $T \cap W$, and consider the dichotomy contained in Lemma \ref{annuli}. If $A$ is parallel to an annulus in $\partial W - {\rm int} N(K)$, then we can decrease $|T \cap W|$ by an isotopy of $T$. We may thus assume every component of $T \cap W$ is an annulus parallel in $W - {\rm int} N(K)$ to $\partial N(K) \cap W$. In which case, $T \cap \partial W$ consists of loops parallel to $K$ in $\partial W$.

Now let $A$ denote any component of $T \cap f(V)$. By the preceding argument, both components of $\partial A$ are parallel to $K$ in $\partial W$. Hence, $A$ is parallel in the handlebody $f(V)$ to the annulus $A'$ on $\partial f(V)$ bounded by $\partial A$. By the minimality of $|T \cap W|$, so $A'$ necessarily contains $K$.

We conclude that $T$ is the union of two annuli, one properly embedded in $f(V)$ and the other properly embedded in $W$ and both parallel to $A'$. It follows that $T$ is necessarily peripheral in $E$, and hence $K$ is atoroidal. $\square$\\

The set $k \cup \alpha \cup \beta$ divides $\partial V$ into a pair of $3$-holed spheres, $P_{1}$ and $P_{2}$. We now define $S_{ij}$ to be equal to $f(P_{i} \cup A_{j})$, for each $i, j \in \{1, 2\}$. Each is a genus $1$ Seifert surface for $K$ and, by Lemma \ref{nontriv}, each is therefore of minimal genus. Reindexing if need be, we may assume $S_{11}$ and $S_{22}$ intersect transversely along $\alpha$ and along $\beta$. Let us abbreviate $S_{ii}$ to $S_{i}$ for both $i \in \{1, 2\}$. Then, $S_{1} \cup S_{2}$ divides $\B{S}^{3}$ into the following three regions:

\begin{itemize}
\item $W$, a hyperbolic $3$-manifold;
\item The solid torus $f(V_{0})$, bounded by $f(A_{1})$ and $f(A_{2})$, and
\item A third region that contains $S_{1} \cap S_{2}$ and that is branched along $S_{1} \cap S_{2}$. In particular, this region is not a $3$-manifold.
\end{itemize}

\noindent None of these regions can give a product region between $S_{1}$ and $S_{2}$. (Recall $f(V_{0})$ does not give a product region between $A_{1}$ and $A_{2}$.) It follows from the contrapositive of Proposition 4.8(2) in [12] that $S_{1}$ and $S_{2}$ intersect essentially. In particular, $d([S_{1}], [S_{2}]) \geq 2$. According to Corollary \ref{atmost2}, the diameter of $\mathcal{G}(K)$ is at most $2$. We conclude that the diameter of $\mathcal{G}(K)$ is exactly $2$. This completes a proof of Proposition \ref{sharp}.

\section{An infinite class of atoroidal knots}

The purpose of this section is to prove Propostion \ref{seventh}, offering a family of atoroidal knots, parameterised by knot genus, each of whose associated graphs has diameter precisely the modulus of the knot Euler characteristic. In particular, their diameters grow linearly with knot genus.

Given any non-negative integer $g$, pick a sequence of integers $a_1, a_2, \ldots, a_{g}$ of length $g$ such that $|a_j|\ge 2$ for every $j$. Let $K$ be the 2-bridge knot whose slope is represented by the continued fraction\[
[2a_1,-2a_2, \ldots, 2a_{2g-1}, -2a_{2g}]
=
\cfrac{1}{2a_1-
\cfrac{1}{2a_2 + \ldots
\cfrac{1}{2a_{2g-1}-
\cfrac{1}{2a_{2g}}}}}
\]
Then, the genus of $K$ is precisely $g$. We show that the diameter of $\mathcal{G}(K)$ is equal to $2g - 1$ by using [12], where the structure of Kakimizu's complex $\mathcal{MS}(K)$ is explicitly described. To recall, let $T$ be a tree, with $n:=2g$ vertices, whose underlying space is homeomorphic to a closed interval, and let $v_1, v_2, \ldots, v_n$ be the vertices of $T$, lying on the interval in this order. For each vertex $v_j$ we associate an unknotted oriented annulus $F(v_j)$ in $S^3$ with $a_j$-right hand full twists. Then, $K$ is equal to the boundary of a surface
obtained by successively plumbing the annuli $F(v_1), F(v_2), \ldots, F(v_n)$, and this surface is a minimal genus Seifert surface for $K$. Moreover, every minimal genus Seifert surface of $K$ is obtained in this way (see [2]).

There are $2^{n-1}$ different ways of successive plumbing, according as $A_{j+1}$ is plumbed to $A_j$ from above or from below with respect to a normal vector field on $A_{j}$. Thus, successive plumbing can be represented by an {\em orientation of $T$}, directing each edge in one of two ways, by the following rule: If $\rho$ is an orientation of $T$, then we plumb $A_{j+1}$ to $A_j$ from above or below according as the edge joining $v_j$ and $v_{j+1}$ has initial point $v_j$ or $v_{j+1}$, respectively, with respect to $\rho$. See Section 2 of [12] for a more detailed account.

We denote by $S(\rho)$ the Seifert surface of $K$ determined by the orientation $\rho$.  The condition that $|a_j|\ge 2$ for every $j$ guarantees the correspondence $\rho \mapsto S(\rho)$ determines a bijection from the set $\OO(T)$ of all orientations of $T$ to the vertex set of $\mathcal{MS}(K)$ (see Theorem 2.3 of [12]).

To describe the structure of $\mathcal{MS}(K)$, we introduce a few definitions. A vertex $v_{j}$ of $T$ is said to be a {\em sink} for the orientation $\rho$ of $T$ if every edge of $T$ incident on $v_{j}$ points towards $v_{j}$. If $v_{j}$ is a sink for $\rho$, then let $v_j(\rho)$ denote the orientation of $T$ obtained from $\rho$ by reversing the orientations of each edge incident on $v_j$. A {\em cycle} in $\OO(T)$ is a sequence

\[
\rho_1 \xrightarrow{v_{j_1}} 
\rho_2 \xrightarrow{v_{j_2}} 
\ldots \xrightarrow{v_{j_{n-1}}}
\rho_{n} \xrightarrow{v_{j_{n}}}
\rho_1,
\] where $(j_1,j_2,\ldots, j_n)$ is a permutation of $\{1,2,\ldots,n\}$ and $\rho_1,\rho_2,\ldots,\rho_n$ are mutually distinct elements of $\OO(T)$ such that $v_{j_{k}}(\rho_k)=\rho_{k+1}$ for every $k$, where our indices are considered modulo $n$. According to Theorem 3.3 of [12], $\mathcal{MS}(K)$ can be described as follows:

\begin{itemize}
\item The vertex set of $\mathcal{MS}(K)$ is identified with $\OO(T)$, and

\item A set of vertices $\{\rho_0,\rho_1,\ldots,\rho_k\}$ spans a $k$-simplex in $\mathcal{MS}(K)$ if and only if it is contained in a cycle of $\mathcal{O}(T)$.
\end{itemize}

Moreover, $\mathcal{MS}(K)$ gives a triangulation of the cube $I^{n-1}$
whose vertices are all the corners of the cube (see Proposition 3.9 of [12]).

We now show that the diameter of $\mathcal{G}(K)$ is equal to $n-1$. Identify $\OO(T)$ with $\{-,+\}^{(n-1)}$ by identifying $\rho\in\OO(T)$ with $(\epsilon_1,\epsilon_2,\ldots, \epsilon_{n-1})$, where $\epsilon_j$ is $+$ or $-$ according as the initial point of the $j$-th edge is $v_j$ or $v_{j+1}$, respectively.

\newtheorem{7.1}[1.1]{Lemma}
\begin{7.1}
For any two elements $\rho$ and $\rho'$ of $\OO(T)$, we have $d(\rho,\rho')\le n-1$, where $d$ denotes the edge-path distance in $\mathcal{G}(K)$.
\label{eighth}
\end{7.1}

\noindent Proof.  We prove the lemma by inducting on $n$. Note that if $n$ is odd, so not of the form $2g(K)$, we may still consider a linear tree $T$ with $n$ vertices and a simplicial complex with vertex set $\OO(T)$. If $n=1$, $\OO(T)$ consists of a single element and the lemma obviously holds.

Let $\rho = (\epsilon_1, \epsilon_2, \ldots, \epsilon_{n-1})$ and $\rho' = (\epsilon_1', \epsilon_2', \ldots, \epsilon_{n-1}')$ be two elements of $\OO(T)$, where $T$ has $n$ vertices. Suppose first that $\epsilon_{n-1}=\epsilon_{n-1}'$. Let $T_0$ be the sub-tree of $T$ obtained by deleting the last edge. By the inductive hypothesis, the distance between $\rho_0:=(\epsilon_1,\epsilon_2,\ldots, \epsilon_{n-2})$ and $\rho_0':=(\epsilon_1',\epsilon_2',\ldots, \epsilon_{n-2}')$ in $\OO(T_0)$ is at most $n-2$.
Since every edge in $\OO(T_0)=\{-,+\}^{(n-2)}$ lifts to an edge in $\{-,+\}^{(n-2)}\times\{\epsilon\} \subset \OO(T)$, where $\epsilon=\epsilon_{n-1}=\epsilon_{n-1}'$, we see $d(\rho,\rho')$ is at most $n-2$.

Suppose next that $\epsilon_{n-1}\ne\epsilon_{n-1}'$. Let $\rho''$ be the element of $\OO(T)$ obtained from $\rho'$ by replacing $\epsilon_{n-1}'$ with $\epsilon_{n-1}$. Then, we have $$d(\rho,\rho') \leq d(\rho,\rho'')+d(\rho'',\rho') \leq n - 2 + 1 = n - 1$$ and this completes a proof of Lemma \ref{eighth}. $\square$\\

\newtheorem{7.2}[1.1]{Lemma}
\begin{7.2}
Let $\rho_-=(-,-,\ldots,-)$ and $\rho_+=(+,+,\ldots,+)$. Then, $d(\rho_-,\rho_+)\ge n-1$.
\label{ninth}
\end{7.2}

\noindent Proof. Let $w(\rho)$ be the number of $+$ entries of $\rho\in \{-,+\}^{(n-1)}$, so that $w(\rho_-)=0$ and $w(\rho_+)=n-1$. The statement of the lemma follows once we prove $|w(\rho)-w(\rho')|\le 1$ for any edge $(\rho,\rho')$ of $\mathcal{G}(K)$. To prove this, observe that if $v_j$ is a sink for $\rho$ then $w(v_j(\rho))$ is equal to $w(\rho)$, $w(\rho)+1$, or $w(\rho)-1$ according as $j\in \{2,3,\ldots,n-1\}$, $j=1$ or $j=n$. Let $(\rho,\rho')$ be an edge of $\mathcal{G}(K)$. Then, $\{\rho,\rho'\}$ is contained in the vertex set of a maximal simplex, of $\mathcal{MS}(K)$, which in turn is the set of all orientations for some cycle, say $$\rho_1 \xrightarrow{v_{j_1}} \rho_2 \xrightarrow{v_{j_2}} \ldots \xrightarrow{v_{j_{n-1}}} \rho_{n} \xrightarrow{v_{j_{n}}} \rho_1.$$ Since every vertex appears in a cycle, the above observation implies the set $\{w(\rho_1),w(\rho_2),\ldots, w(\rho_n)\}$ consists of two successive integers. In particular $|w(\rho)-w(\rho')|\le 1$, and this completes a proof of Lemma \ref{ninth}. $\square$\\

By Lemma \ref{eighth} and Lemma \ref{ninth} we see that the diameter of $\mathcal{G}(K)$ is equal to $n-1$, thus completing a proof of Proposition \ref{seventh}.\\

\begin{center}
R\begin{small}ESEARCH \end{small}U\begin{small}PDATE\end{small}\\
\end{center}

Roberto Pelayo's thesis [11] became publicly available from April 2007. The upper bound given in Theorem 10.1 of [11] is quadratic in knot genus though is not computed. The argument found therein is based on minimal surface theory, and is quite different from that given here in Section 3.

In August 2007, Jennifer Schultens [14] gave an elegant proof of the simple connectivity of Kakimizu's complex, using PL-minimal surface theory. In fact her argument can be extended to prove that the second homotopy group of Kakimizu's complex is also trivial.

The Kakimizu conjecture remains open.


\section*{Acknowledgements}

The authors wish to thank Osamu Kakimizu, for two particularly illuminating conversations, and to thank Tom Fleming, David Futer, Fran\c cois Gu\'eritaud, Mikami Hirasawa, Kazuhiro Ichihara, Toshio Saito, and Yukihiro Tsutsumi for many interesting and helpful conversations. This project began life while the second author was visiting Osaka University. The second author wishes to thank Osaka University for its warm hospitality. 

The first author is supported by a Japan Society for the Promotion of Science grant, number 18340018. The second author was partially supported by a short-term Japan Society for the Promotion of Science post-doctoral fellowship, number PE05043, and is supported by a long-term Japan Society for the Promotion of Science post-doctoral fellowship, number P06034. Both authors wish to thank the JSPS for its financial support.



\begin{tabular}{l l l}
& & Makoto Sakuma\\
& & Department of Mathematics\\
& & Graduate School of Science\\
& & Hiroshima University\\
& & 1-3-1 Kagamiyama\\
& & Higashi-Hiroshima\\
& & 739-8526 Japan\\
& & e-mail: sakuma@math.sci.hiroshima-u.ac.jp\\
& & URL: http://www.math.sci.hiroshima-u.ac.jp/$\sim$sakuma\\
& & \\
& & Kenneth J. Shackleton\\
& & Department of Mathematical\\
& & \indent and Computing Sciences\\
& & Tokyo Institute of Technology\\
& & 2-12-1 O-okayama\\
& & Meguro-ku\\
& & Tokyo\\
& & 152-8552 Japan\\
& & e-mail: shackleton.k.aa@m.titech.ac.jp\\
& & e-mail: kjs2006@alumni.soton.ac.uk\\
& & URL: http://www.is.titech.ac.jp/$\sim$kjshack5\\
\end{tabular}


\begin{thebibliography}{99}

\bibitem{F} S. R. Fenley: \textit{Quasi-Fuchsian Seifert surfaces}, Math. Z. \textbf{228} no. 2 (1998), 221--227.

\bibitem{HaThu} A. E. Hatcher and W. P. Thurston: \textit{Incompressible surfaces in $2$-bridge knot complements}, Invent. Math. \textbf{79} (1985), 225--246.

\bibitem{HiSa} M. Hirasawa and M. Sakuma: \textit{Minimal genus Seifert surfaces for alternating links}, Proceedings of Knots 96 (1997), 383--394.


\bibitem{K0} O. Kakimizu: Talk at the meeting ``Knot theory and related topics'' held in Osaka, 1989.


\bibitem{K2} O. Kakimizu: \textit{Finding disjoint incompressible spanning surfaces for a link}, Hiroshima Math. J. \textbf{22} (1992), 225--236.

\bibitem{K3} O. Kakimizu: \textit{Classification of the incompressible spanning surfaces for prime knots of 10 or less crossings}, Hiroshima Math. J. \textbf{35} no. 1 (2005), 47--92.

\bibitem{Ko1} T. Kobayashi: \textit{Casson-Gordon's rectangle condition of Heegaard diagrams and incompressible tori in $3$-manifolds}, Osaka J. Math. \textbf{25} no. 3 (1988), 553--573.

\bibitem{Ko2} T. Kobayashi: \textit{Uniqueness of minimal genus Seifert surfaces for links}, Topology and its Appl. \textbf{33} (1989) 265--279.

\bibitem{O} U. Oertel: \textit{On the existence of infinitely many essential surfaces of bounded genus}, Pacific J. Math. \textbf{202} no. 2 (2002), 449--458.

\bibitem{P1} R. C. Pelayo: e-mail communication, October 2006.

\bibitem{P2} R. C. Pelayo: \textit{Diameter bounds on the complex of minimal genus Seifert surfaces for hyperbolic Knots}, PhD thesis, April 2007. Available online from the Caltech Library System at http://resolver.caltech.edu/CaltechETD:etd-06042007-015951.

\bibitem{Sa} M. Sakuma: \textit{Minimal genus Seifert surfaces for special aborescent links}, Osaka J. Math. \textbf{31} (1994), 861--905.

\bibitem{ScTho} M. G. Scharlemann and A. A. Thompson: \textit{Finding disjoint Seifert surfaces}, Bull. London Math. Soc. \textbf{20} (1988) 61--64.

\bibitem{Schu} J. C. Schultens: \textit{Covering spaces and the Kakimizu complex}, arXiv:0707.3926v3. 

\bibitem{Thur1} W. P. Thurston: \textit{A norm for the homology of $3$-manifolds}, Mem. Amer. Math. Soc. \textbf{59} no. 339 (1986) i--vi and 99--130.

\bibitem{Thur2} W. P. Thurston: \textit{Three-dimensional geometry and topology}, Princeton University Press, 1997.

\bibitem{Tsu1} Y. Tsutsumi: \textit{Universal bounds for genus one Seifert surfaces for hyperbolic knots and surgeries with non-trivial JSJT-decompositions}, in Proceedings of the Winter Workshop of Topology/Workshop of Topology and Computer (Sendai, 2002/Nara, 2001), Interdisciplinary Information Sciences \textbf{9} no. 1 (2003), 53--60.

\bibitem{Tsu2} Y. Tsutsumi: \textit{Hyperbolic knots with a large number of disjoint minimal genus Seifert surfaces}, preprint.

\bibitem{Wa} F. Waldhausen: \textit{On irreducible $3$-manifolds which are sufficiently large}, Ann. of Math. \textbf{87} (1968) 56--88.

\bibitem{Wi} R. T. Wilson: \textit{Knots with infinitely many incompressible Seifert surfaces}, arXiv:math.GT/0604001.\\~\\

\end{thebibliography}
\end{document}